\documentclass{amsart}
\usepackage{amsfonts}
\usepackage{latexsym}
\usepackage{amssymb}
\usepackage{amsmath}

\newcommand{\R}{\mathbb R}

\newtheorem{thm}{Theorem}[section]

\newtheorem{proposition}{Proposition}[section]
\theoremstyle{remark}
\newtheorem*{rmk}{Remark}

\begin{document}


\title{On an extension of the Blaschke-Santal\'{o} inequality}

\author{David Alonso-Guti\'{e}rrez*}

\thanks{*Supported by FPU Scholarship from MEC (Spain), MCYT
Grants(Spain)  MTM2007-61446, DGA E-64 and by Marie Curie RTN
CT-2004-511953}

\email{498220@celes.unizar.es}

\address{Universidad de Zaragoza}

\date{October 2007}

\begin{abstract}
Let $K$ be a convex body and $K^\circ$ its polar body. Call $\phi(K)=\frac{1}{|K||K^\circ|}\int_K\int_{K^\circ}\langle x,y\rangle^2 dxdy$. It is conjectured that $\phi(K)$ is maximum when $K$ is the euclidean ball. In particular this statement implies the Blaschke-Santal\'{o} inequality. We verify this conjecture when $K$ is restricted to be a $p$--ball.
\end{abstract}

\maketitle

\section{Introduction and notation}

\bigskip
A convex body $K\subset \R^n$ is a compact convex set with non-empty interior. For every convex body, its polar set is defined
$$
K^\circ=\{x\in\R^n\,:\,\langle x,y\rangle\leq 1 \textrm{ for all }y\in K\}
$$
where $\langle \cdot,\cdot\rangle$ denotes the standard scalar product in $\R^n$. Note that if $0\in \textrm{int}K$ then $K^\circ$ is a convex body.

\bigskip
For $p\in[1,\infty]$, let us denote by $B_p^n$ the unit ball of the $p$--norm. It is:
$$
B_p^n=\left\{x\in\R^n\,:\,\sum_{i=1}^n|x_i|^p\leq 1\right\}\hskip 2cm B_\infty^n=\{x\in\R^n\,:\,\max|x_i|\leq 1\}.
$$
It is well known that the polar body of $B_p^n$ is $B_q^n$ where $q$ is the dual exponent of $p$ $(\frac{1}{p}+\frac{1}{q}=1)$. Along this paper $q$ will always denote the dual exponent of $p$.

\bigskip
Given two symmetric convex bodies $A\subset\R^n$, $B\subset\R^m$, for any $p\in[1,\infty]$ they define a symmetric convex body $A\times_p B\subset\R^{n+m}$ which is the unit ball of the norm given by
$$
\Vert(x_1,x_2)\Vert_{A\times_p B}^p=\Vert x_1\Vert_A^p+\Vert x_2\Vert_B^p\hskip 1.8cm\Vert(x_1,x_2)\Vert_{A\times_\infty B}=\max\{\Vert x_1\Vert_A,\Vert x_2\Vert_B\}.
$$
Note that the polar body of $A\times_p B$ is $A^\circ\times_q B^\circ$ and $B_p^n=B_p^{n-1}\times_p[-1,1]$.

\bigskip
A convex body $K$ is said to be in isotropic position if it has volume 1 and satisfies the following two conditions:
\begin{itemize}
\item$\int_Kxdx=0 \textrm{ (center of mass at 0)}$
\item$\int_K\langle x,\theta\rangle^2dx=L_K^2\quad \forall \theta\in S^{n-1}$
\end{itemize}
where $L_K$ is a constant independent of $\theta$, which is called the isotropy constant of $K$.

\bigskip
We will use the notation $\widetilde{K}$ for $|K|^{-\frac{1}{n}}K$.

\bigskip
Given a centrally symmetric convex body $K$, we call
$$
\phi(K)=\frac{1}{|K||K^\circ|}\int_K\int_{K^\circ}\langle x,y\rangle^2 dxdy.
$$
Note that $\phi(K)=\phi(TK)$ for all $T\in GL(n)$. It is conjectured in \cite{KUP} that $\phi(K)$ is maximized by ellipsoids. It is, for every symmetric convex body $K\subset\R^n$
\begin{eqnarray}\label{conjecture}
\phi(K)\leq\phi(B_2^n)=\frac{n}{(n+2)^2}.
\end{eqnarray}

\bigskip
\begin{rmk}
We can also define the functional $\phi$ when $K$ is not symmetric. When $K$ is a regular simplex with its center of mass at the origin, it is easy to compute that $\phi(K)=\phi(B_2^n)$.
\end{rmk}

\bigskip
The Blaschke-Santal\'{o} inequality \cite{SA} says that for every symmetric convex body $K$
$$
|K||K^\circ|\leq |B_2^n|^2.
$$
The conjecture (\ref{conjecture}) is stronger than the Blaschke-Santal\'{o} inequality since
$$
\frac{n|K|^\frac{2}{n}|K^\circ|^\frac{2}{n}}{(n+2)^2|B_2^n|^\frac{4}{n}}\leq\frac{1}{|K||K^\circ|}\int_K\int_{K^\circ}\langle x,y\rangle^2 dxdy.
$$
This fact is a consequence of Lemma 6 in \cite{BALL}. In \cite{BALL2}, Ball proved that for 1-unconditional bodies
$$
\int_K\int_{K^\circ}\langle x,y\rangle^2 dxdy\leq \frac{n|B_2^n|^2}{(n+2)^2}
$$
and suggested that this inequality might be true for every convex body. This assertion is slightly weaker than the conjecture in \cite{KUP}, which is not known to be true even for 1-unconditional bodies. In section 2 we are going to prove that the conjecture is true if we restrict $K$ to be a $p$--ball, for some $p\geq 1$. We state this as a theorem:

\begin{thm}\label{teorema}
Among the $p$--balls, the functional $\phi$ is maximized for the euclidean ball .
$$
\max_{p\in[1,\infty]}\phi(B_p^n)=\phi(B_2^n)=\frac{n}{(n+2)^2}.
$$
\end{thm}

\bigskip
The conjecture (\ref{conjecture}) is also stronger than the hyperplane conjecture, which says that there exists an absolute constant $C$ such that for every isotropic convex body $L_K<C$. It can be proved that $\phi(K)$ is bounded from below by $\frac{c_1}{n}$, where $c_1$ is an absolute constant. If there exists an absolute constant $c_2$ such that $\phi(K)\leq\frac{c_2}{n}$, then the hyperplane conjecture would be true, since
$$
nL_K^2L_{K^\circ}^2\leq \frac{\phi(K)}{|K|^\frac{2}{n}|K^\circ|^\frac{2}{n}}\leq cn^2\phi(K)
$$
where $c$ is an absolute constant.

\bigskip

In case that $\widetilde{K}$ and $\widetilde{K^\circ}$ are both isotropic then $\phi(K)=n|K|^\frac{2}{n}|K^\circ|^\frac{2}{n}L_K^2L_{K^\circ}^2$ and the conjecture $\phi(K)\sim\frac{1}{n}$ is equivalent to the hyperplane conjecture. This is the case of 1-symmetric bodies, for which the hyperplane conjecture is known to be true (A convex body is 1-symmetric if it is invariant under reflections in the coordinate hyperplanes and under permutations of the coordinates).

\bigskip

We will say that a symmetric convex body $K\subset\R^n$ is a revolution body if there exists $\theta\in S^{n-1}$ and a concave function $r(t)$ such that for every $t\in[-h_K(\theta),h_K(\theta)]$ $K\cap(t\theta +\theta^\perp)=r(t)B_2^{n-1}$, where $h_K(\theta)$ is the support function of $K$:
$$
h_K(\theta)=\max\{\langle x,\theta\rangle\,:\, x\in K\}.
$$

\bigskip
In section 3 we will prove that there exists an absolute constant $C$ such that whenever $K$ is a symmetric convex body of revolution, $\phi(K)\leq\frac{C}{n}$.

\bigskip
Along this paper, $\psi$ will always denote the logarithmic derivative of the Gamma function. We will make use of the following identity on the derivatives of $\psi$, known as polygamma functions:
$$
\psi^{(n)}(x)=(-1)^{n+1}\int_0^\infty t^n\frac{e^{-xt}}{1-e^{-t}}dt.
$$

\bigskip
The letters $C,c_1,c_2,\dots$ will always denote absolute constants which do not depend on the dimension.

\bigskip
\section{The $p$--balls}

\bigskip
In this section we are going to prove theorem \ref{teorema}. We will obtain it as a consequence of the following
\begin{thm}
For every $A\subset\R^n$, $B\subset\R^m$, $p\in[1,\infty]$
$$
\phi(A\times_pB)=f(n,n+m,p)\phi(A)+f(m,n+m,p)\phi(B)
$$
where
$$
f(y_1,y_2,p)=\left\{\begin{array}{lr}\frac{(y_1+2)^2y_2^2\Gamma\left(\frac{y_1+2}{p}\right)\Gamma\left(\frac{y_1+2}{q}\right)\Gamma\left(\frac{y_2}{p}\right)\Gamma\left(\frac{y_2}{q}\right)}
{y_1^2(y_2+2)^2\Gamma\left(\frac{y_2+2}{p}\right)\Gamma\left(\frac{y_2+2}{q}\right)\Gamma\left(\frac{y_1}{p}\right)\Gamma\left(\frac{y_1}{q}\right)} &p,\neq 1,\infty\\ \\ \frac{(y_1+2)y_2\Gamma\left(y_1+2\right)\Gamma\left(y_2\right)}
{y_1(y_2+2)\Gamma\left(y_2+2\right)\Gamma\left(y_1\right)} &p=1,\infty\end{array} \right.
$$
attains its maximum when $p=2$, for every $0<y_1<y_2$.
\end{thm}

\bigskip
\begin{proof}
First of all we are going to prove that for every fixed $0<y_1<y_2$, the function defined on $[0,1]$ like $f_1(x)=f(y_1,y_2,\frac{1}{x})$ attains its maximum in $x=\frac{1}{2}$. It is easy to check that $f_1(0)=f_1(1)<f_1(\frac{1}{2})$. $f_1$ has got a maximum in $x=\frac{1}{2}$ if and only if $\log f_1$ has got a maximum in $x=\frac{1}{2}$.

\bigskip
Since $f_1(x)=f_1(1-x)$, it is enough to prove that $\log f_1$ is increasing in $(0,\frac{1}{2})$. Now, if we call
$$
F(x,y)=(y+2)[\psi((y+2)x)-\psi((y+2)(1-x))]-y[\psi(yx)-\psi(y(1-x))]
$$
we have that
\begin{eqnarray*}
(\log f_1)'(x)&=&(y_1+2)[\psi((y_1+2)x)-\psi((y_1+2)(1-x))]-y_1[\psi(y_1x)-\psi(y_1(1-x))]\\
& &-(y_2+2)[\psi((y_2+2)x)-\psi((y_2+2)(1-x))]+y_2[\psi(y_2x)-\psi(y_2(1-x))]\\
&=&F(x,y_1)-F(x,y_2).
\end{eqnarray*}
So it is enough to prove that for every fixed $x\in(0,\frac{1}{2})$, $F(x,y)$ is decreasing in $y\in(0,\infty)$. Hence we compute
\begin{eqnarray*}
\frac{\partial F}{\partial y}(x,y)&=&\psi((y+2)x)-\psi((y+2)(1-x))-\psi(yx)+\psi(y(1-x))\\
& &+(y+2)x\psi'((y+2)x)-(y+2)(1-x)\psi'((y+2)(1-x))\\
& &-yx\psi'(yx)+y(1-x)\psi'(y(1-x)).\\
\end{eqnarray*}

\bigskip
We call this last quantity $G(x,y)$ and we will see that $G(x,y)<0$ if $x\in(0,\frac{1}{2})$ and $G(x,y)>0$ if $x\in(\frac{1}{2},1)$. Notice that $G(\frac{1}{2},y)=0$, so we just need to check that for every fixed $y>0$, $G(x,y)$ is increasing in $x$. Computing its derivative we obtain
\begin{eqnarray*}
\frac{\partial G}{\partial x}(x,y)&=&2(y+2)[\psi'((y+2)x)+\psi'((y+2)(1-x))]\\
& &+(y+2)^2[x\psi''((y+2)x)+(1-x)\psi''((y+2)(1-x))]\\
& &-2y[\psi'(y(1-x))+\psi'(yx)]-y^2[x\psi''(yx)+(1-x)\psi''(y(1-x))]\\
&=&H(x,y+2)-H(x,y).
\end{eqnarray*}
where we have called $H(x,y)$ the following function
$$
H(x,y)=2y[\psi'(yx)+\psi'(y(1-x))]+y^2[x\psi''(yx)+(1-x)\psi''(y(1-x))].
$$

\bigskip
Hence, if for every fixed $x\in(0,1)$ $H(x,y)$ is increasing in $y$, then so it is $G(x,y)$ in $x$ for fixed $y$ and the theorem is proved. In order to prove this, we need the following result concerning the $\psi$ function whose proof can be found in \cite{ALZ}. We will write it here for the sake of completeness:

\bigskip
\begin{proposition}\label{conv}
The function $f(x)=x^2\psi'(x)$ is convex in the interval $(0,\infty)$.
\end{proposition}
\begin{proof}
The second derivative of $f$ is
$$
f''(x)=2\psi'(x)+4x\psi''(x)+x^2\psi'''(x).
$$

\bigskip
Using the integral representation of the derivatives of $\psi$ this is equal to
\begin{eqnarray*}
f''(x)&=&\int_0^\infty\frac{e^{-xt}}{1-e^{-t}}(2t-4xt^2+x^2t^3)dt\\
&=&\int_0^\infty\frac{t}{1-e^{-t}}\frac{d^2}{dt^2}(t^2e^{-xt})dt\\
&=&\int_0^\infty \frac{d^2}{dt^2}\left(\frac{t}{1-e^{-t}}\right)t^2e^{-xt}dt
\end{eqnarray*}
which is positive since the function $\frac{t}{1-e^{-t}}$ is convex in the interval $(0,\infty)$.
\end{proof}

\bigskip
Now, for every $x\in(0,1)$, $y>0$ we have that
\begin{eqnarray*}
\frac{\partial H}{\partial y}(x,y)&=&2\psi'(yx)+4yx\psi''(yx)+y^2x^2\psi'''(yx)\\
&+&2\psi'(y(1-x))+4y(1-x)\psi''(y(1-x))+y^2(1-x)^2\psi'''(y(1-x))>0
\end{eqnarray*}
as a consequence of proposition \ref{conv} and this proves that $f(y_1,y_2,p)\leq f(y_1, y_2,2)$ when $0<y_1<y_2$.

\bigskip
Let us prove now that
$$
\phi(A\times_pB)=f(n,n+m,p)\phi(A)+f(m,n+m,p)\phi(B).
$$
Assume that $p\neq1,\infty$. We compute the volume of $A\times_p B$:
\begin{eqnarray*}
|A\times_p B|&=&\int_A(1-\Vert x_1\Vert_A^p)^{\frac{m}{p}}|B|dx_1=\int_A\int_{\Vert x_1\Vert_A^p}^1\frac{m}{p}(1-t)^{\frac{m}{p}-1}dt|B|dx_1\cr
&=&\int_0^1\int_{ t^\frac{1}{p}A}\frac{m}{p}(1-t)^{\frac{m}{p}-1}|B|dx_1dt=\frac{m}{p}|A||B|\beta\left(\frac{m}{p}+1,\frac{n}{p}\right)\cr
&=&\frac{nm}{p(n+m)}|A||B|\beta\left(\frac{m}{p},\frac{n}{p}\right).
\end{eqnarray*}

\bigskip
Since $(A\times_p B)^\circ =A^\circ\times_q B^\circ$, we have that
$$
|(A\times_p B)^\circ|=\frac{nm}{q(n+m)}|A^\circ||B^\circ|\beta\left(\frac{m}{q},\frac{n}{q}\right).
$$

\bigskip
From the symmetry of $A$ and $B$ we obtain that
$$
\int_K\int_{K^\circ}\langle (x_1,x_2),(y_1,y_2)\rangle^2 dydx=\int_K\int_{K^\circ}\langle x_1,y_1\rangle^2 dydx +\int_K\int_{K^\circ}\langle x_2,y_2\rangle^2 dydx
$$
where we have called $K=A\times_p B$.

\bigskip
Let us compute these integrals:
\begin{eqnarray*}
&&\int_K\int_{K^\circ}\langle x_1,y_1\rangle^2dydx\cr
&=&\int_A\int_{A^\circ}\langle x_1, y_1\rangle^2(1-\Vert x_1\Vert_A^p)^\frac{m}{p}(1-\Vert y_1\Vert_{A^\circ}^q)^\frac{m}{q}|B||B^\circ|dy_1dx_1\cr
&=&|B||B^\circ|\int_A\int_{A^\circ}\langle x_1,y_1\rangle^2\int_{\Vert x_1\Vert_A^p}^1\frac{m}{p}(1-t)^{\frac{m}{p}-1}dt\int_{\Vert y_1\Vert_{A^\circ}^q}^1\frac{m}{q}(1-s)^{\frac{m}{q}-1}dsdy_1dx_1\cr
&=&|B||B^\circ|\frac{m^2}{pq}\int_0^1\int_0^1\int_{t^\frac{1}{p}A}\int_{s^\frac{1}{q}A^\circ}\langle x_1,y_1\rangle^2(1-t)^{\frac{m}{p}-1}(1-s)^{\frac{m}{q}-1}dy_1dx_1dsdt\cr
&=&|B||B^\circ|\frac{m^2}{pq}\beta\left(\frac{m}{p},\frac{n+2}{p}+1\right)\beta\left(\frac{m}{q},\frac{n+2}{q}+1\right)\int_A\int_{A^\circ}\langle x_1,y_1\rangle^2dy_1dx_1\cr
&=&|B||B^\circ|\frac{m^2(n+2)^2}{pq(m+n+2)^2}\beta\left(\frac{m}{p},\frac{n+2}{p}\right)\beta\left(\frac{m}{q},\frac{n+2}{q}\right)\int_A\int_{A^\circ}\langle x_1,y_1\rangle^2dy_1dx_1\cr
\end{eqnarray*}
and in the same way
\begin{eqnarray*}
&&\int_K\int_{K^\circ}\langle x_2,y_2\rangle^2dydx=\cr
&=&|A||A^\circ|\frac{n^2(m+2)^2}{pq(m+n+2)^2}\beta\left(\frac{n}{p},\frac{m+2}{p}\right)\beta\left(\frac{n}{q},\frac{m+2}{q}\right)\int_B\int_{B^\circ}\langle x_1,y_1\rangle^2dy_1dx_1.
\end{eqnarray*}

\bigskip
Now from the definition of $\phi$ and the identity $\beta(x,y)=\frac{\Gamma(x)\Gamma(y)}{\Gamma(x+y)}$ we obtain the result. When $p=1,\infty$ the theorem is proved in the same way.
\end{proof}

\bigskip
\section{Revolution bodies}

\bigskip
In this section we are going to prove the following:
\begin{thm}
There exists an absolute constant $C$ such that for every symmetric convex body of revolution $K$, $\phi(K)< \frac{C}{n}$.
\end{thm}

This is not a new result since A. Giannopoulos proved it in his PhD thesis but it was left unpublished. I would like to thank him for allowing me to add this result to this paper.

\bigskip

\begin{proof}
Since $\phi(TK)=\phi(K)$ for every $T\in GL(n)$, we can assume that
$$K=\{\bar{x}=(t,x)\in\R^n\,:\,t\in[-1,1]\,,\,|x|\leq r_1(t)\}$$
where $r_1(t)$ is a concave function such that $r_1(0)=1$.

\bigskip
Then, $K^\circ$ is another revolution body
\begin{eqnarray*}
K^\circ&=&\{\bar{y}=(s,y)\in\R^n\,:\,ts+r_1(t)|y|\leq 1\,,\, \forall t\in[-1,1]\}\cr
&=&\{\bar{y}=(s,y)\in\R^n\,:\,s\in[-1,1]\,,\,|y|\leq r_2(s)\}
\end{eqnarray*}
where $r_2(s)$ is a concave function such that $r_2(0)=1$

\bigskip
Let us now compute $\phi(K)$:
\begin{eqnarray*}
\phi(K)&=&\frac{1}{|K||K^\circ|}\int_K\int_{K^\circ}(ts+\langle x,y\rangle)^2d\bar{y}d\bar{x}\cr
&=&\frac{1}{|K||K^\circ|}\int_K\int_{K^\circ}t^2s^2+\langle x,y\rangle^2d\bar{y}d\bar{x}\cr
&=&\frac{1}{|K||K^\circ|}\int_K\int_{K^\circ}t^2s^2d\bar{y}d\bar{x}\cr
&&+\frac{1}{|K||K^\circ|}\int_{-1}^1\int_{-1}^1\int_{r_1(t)B_2^{n-1}}\int_{r_2(t)B_2^{n-1}}\langle x,y\rangle^2dydxdsdt\cr
&=&|K|^\frac{2}{n}|K^\circ|^\frac{2}{n}\int_{\widetilde{K}}t^2d\bar{x}\int_{\widetilde{K^\circ}}s^2d\bar{y}\cr
&&+\frac{\int_{-1}^1r_1(t)^{n+1}dt\int_{-1}^1r_2(s)^{n+1}ds}{\int_{-1}^1r_1(t)^{n-1}dt\int_{-1}^1r_2(s)^{n-1}ds}\phi(B_2^{n-1})
\end{eqnarray*}

Since $\max\{r_1(t)\,,\,t\in[-1,1]\}=r_1(0)=1$ and $\max\{r_2(s)\,,\,s\in[-1,1]\}=r_2(0)=1$, for every $t,s\in[-1,1]$ we have that
\begin{itemize}
\item$r_1(t)^{n+1}\leq r_1(t)^{n-1}$
\item$r_2(s)^{n+1}\leq r_2(s)^{n-1}$
\end{itemize}
and hence the second summand is bounded by $\phi(B_2^{n-1})=\frac{n-1}{(n+1)^2}$.

\bigskip
To bound the first summand we will use the following well known result by Hensley\cite{H}:

\bigskip
\textit{``There exist absolute constants $c_1$, $c_2$ such that for every symmetric convex body $K\subset\R^n$ with volume $1$ and for every $\theta\in S^{n-1}$"}
$$
\frac{c_1}{|K\cap\theta^\perp|}\leq\left(\int_K\langle x,\theta\rangle^2dx\right)^\frac{1}{2}\leq\frac{c_2}{|K\cap\theta^\perp|}.
$$

\bigskip
Hence
\begin{itemize}
\item$\int_{\widetilde{K}}t^2d\bar{x}\leq \frac{c|K|^{2\frac{n-1}{n}}}{|K\cap e_1^\perp|^2}=\frac{c|K|^{2-\frac{2}{n}}}{|B_2^{n-1}|^2}$,
\item$\int_{\widetilde{K^\circ}}s^2d\bar{y}\leq \frac{c|K^\circ|^{2\frac{n-1}{n}}}{|K^\circ\cap e_1^\perp|^2}=\frac{c|K^\circ|^{2-\frac{2}{n}}}{|B_2^{n-1}|^2}$.
\end{itemize}

\bigskip
So, by Blaschcke-Santal\'{o} inequality, the first summand is bounded by
$$
\frac{c|K|^2|K^\circ|^2}{|B_2^{n-1}|^4}\leq \frac{c|B_2^n|^4}{|B_2^{n-1}|^4}.
$$

\bigskip
Now, using the fact that $|B_2^n|=\frac{\pi^\frac{n}{2}}{\Gamma\left(1+\frac{n}{2}\right)}$ and Stirling's formula, we obtain that the first summand is bounded by $\frac{c}{n^2}$ and hence the theorem is proved.
\end{proof}
\bigskip
{\centerline {ACKNOWLEDGEMENTS}}

\bigskip
I would  like to thank professor Apostolos Giannopoulos for encouraging me to work on this problem and Greg Kuperberg for helping me to present the result in a more general way.



\begin{thebibliography}{99?}
\bibitem{ALZ} Alzer H.
\textit{Inequalities for the gamma function}
Proceedings of the AMS {\bf 128}, no. 1, pp.141--147

\bibitem{BALL} Ball K. M.
\textit{Logarithmically concave functions and sections of convex sets.}
Studia Math.{\bf 88} (1988), no. 1, pp. 69--84.

\bibitem{BALL2} Ball K. M.
\textit{Some remarks on the geometry of convex sets.}
Geometric aspects of functional analysis (1986/87), pp. 224--231,
Lecture Notes in Math. {\bf 1317}, Springer, Berlin, 1988.


\bibitem{H} Hensley D.
\textit{Slicing convex bodies, bounds of slice area in terms of
the body's covariance}, {Proc. Amer. Math. Soc.} {\bf 79} (1980),
pp.\,619-625.

\bibitem{KUP} Kuperberg G.
\textit{From the Mahler conjecture to Gauss linking integrals.}
Geometric Aspects of Functional Analysis. To appear.

\bibitem{SA} Santal\'{o}, L.
\textit{Un invariante af\'{\i}n para los cuerpos convexos del espacio de $n$ dimensiones.}
Portugal Math. {\bf 8} (1949), pp.\,155-161

\end{thebibliography}
\end{document}